\documentstyle[11pt]{article}
\def\ba{\begin{eqnarray}}
\def\ea{\end{eqnarray}}
\def\lb{\label}
\def\be{\begin{equation}}
\def\ee{\end{equation}}

\def\id{{\bf I}}
\def\C{\!\bf{C}\!\!\!\rule{0.2mm}{2.6mm}\;\,}
\def\u{u_{|1 2 \dots p{\cal i}}}
\def\vv{v^{{\cal h}1 2 \dots p|}}

\begin{document}
\title{Hecke Symmetries and Characteristic Relations
on Reflection Equation Algebras}
\author{D.I. Gurevich\\
{\small \it ISTV, Universit\'e de Valenciennes, 59304 Valenciennes, France}\\
P.N. Pyatov\\
{\small \it Bogoliubov Laboratory of Theoretical Physics, JINR, 141980  Dubna,
Moscow region, Russia} \\
P.A. Saponov\\
{\small\it Theory Department, IHEP, 142284 Protvino, Moscow region, Russia}}
\date{}
\maketitle
\begin{abstract}
We discuss how properties of Hecke symmetry (i.e., Hecke type
${\cal R}$-matrix) influence the algebraic structure of the corresponding
Reflection Equation (RE) algebra.
Analogues of the Newton relations and Cayley-Hamilton theorem
for the matrix of generators of the RE algebra
related to a finite rank even Hecke symmetry
are derived.
\end{abstract}
\section{Introduction}

There exist (at least) two matrix type algebras connected to the quantum
groups. One of them is well known since the creation of the quantum groups
theory. It is dual to a quantum linear group and is defined by the famous
"RTT"  relation (\ref{RTT}). We restrict ourselves here to an
${\cal R}$-matrix
of Hecke type (called the Hecke symmetry below). In this case the
algebra under question is a flat deformation of the corresponding
classical object.

Another "matrix algebra" is related to the so-called reflection
equation (\ref{RLRL}) (RE algebra in what follows).
It was introduced by I.Cherednik
in the context of factorizable scattering on
a half-line \cite{Chered} and then found a
number of different applications (see e.g. \cite{Mog}).
%(for R-matrices with a parameter) and is used in the construction of
%the so-called braided groups by S. Majid.
Both these algebras play an
important role in the integrable systems theory and, moreover there
exists a "transmutation procedure" \cite{Mog2} converting one
of them to another.

For the latter algebra there exists an analogue of Newton relations
between two natural families of invariant "functions".  The principle
aim of these notes is to present these relations and demonstrate a
quantum analogue of Cayley-Hamilton theorem in the RE algebra.

In \cite{PS} and \cite{NT} the RE algebras related to the standard
quantum linear group were considered from this point of view.
Nevertheless, there exists a great deal of non-deformational Hecke
symmetries constructed in \cite{Gur}.  If $R: V^{\otimes 2} :\to
V^{\otimes 2}$ is such a symmetry acting in some linear space $V$ it
naturally generates  "symmetric" and "skew-symmetric" algebras of the
space $V$. Then the dimensions of homogeneous components of these
algebras can be completely different from those of the  classical (or
super-) case. We show here that some versions of the Newton and
Cayley-Hamilton relations are valid for any even (this means that the
skew-symmetric algebra is finite dimensional) and closed (see  below)
Hecke symmetry.

Let us note that for the involutive ${\cal R}$-matrix ($R^2=1$) the
RE algebras were considered under the name of monoidal groups in
\cite{Gur}, where a version of the Newton relations was established for
the simplest (i.e., of rank 2, see below) ${\cal R}$-matrix.

The paper is organized as follows. Section 2 contains the definitions
and a survey of known facts about Hecke symmetries and RE algebras.  In
Section 3 we construct iterative relations connecting two different
sets of the central elements of RE algebras.  The characteristic
identities for the matrix of generators of RE algebras are also derived
here.
These results may be viewed as a generalization of the
classical Newton relations and Cayley-Hamilton theorem (see e.g.
\cite{Lank}) to the case of matrices with non-commuting entries.

\section{Preliminaries}
\setcounter{equation}0

\noindent
{\bf Hecke symmetries and related algebras.}

First, let us recall several basic definitions.
The matrix
$R_{12}
\equiv R_{i_1 i_2}^{\, j_1 j_2} \in Mat(N^2,\C)$,
$i_1,i_2,j_1,j_2 = 1, \dots ,N$,
is called the {\em Hecke symmetry}
if it satisfies the
Yang-Baxter equation and  Hecke condition, respectively,\footnote{
Here and in what follows the standard conventions of \cite{FRT}
are used in denoting matrix indices.
}
\ba
\lb{YB}
 R_{12}  R_{23} R_{12} &=&
 R_{23}  R_{12} R_{23} \; , \\
\lb{H}
R^2 = \id + \lambda  R \; , &&  \lambda = q - 1/q \; .
\ea
Here $q \in \C$ is an arbitrary parameter.
We will assume later on that $q$ is not a root of unity
and, hence, for any natural $p$ the corresponding $q$-number
$p_q\equiv (q^p - q^{-p})/(q-q^{-1})$ is different from zero.

With any Hecke symmetry one associates in a canonical way
the following two unital associative matrix type algebras.

\noindent
{\bf A)}. The algebra ${\cal T}(R)$ \cite{FRT}
is generated by the set of $N^2$
generators $T_i^{\, j}$ subject to the relation
\be
\lb{RTT}
R_{12} T_1 T_2 = T_1 T_2 R_{12} \; .
\ee
It is the Hopf algebra with the comultiplication
\be
\lb{comult}
\triangle T_i^{\, j} = T_i^{\, k} \otimes T_k^{\, j} \; .
\ee
With some additional assumptions on the matrix $R$
(see \cite{Gur} and below)
one can define also the antipodal mapping ${\cal S}(\cdot)$ such that
$$
{\cal S}(T_i^{\, k}) T_k^{\, j} = T_i^{\, k} {\cal S}(T_k^{\, j}) =
{\delta}_i^{\, j} \; .
$$
In case of $R$ being the ${\cal R}$-matrix corresponding to
linear quantum group, the algebra ${\cal T}(R)$ becomes
the quantization of the algebra of functions on the
linear group \cite{FRT}.

\noindent
{\bf B)}. The {\em Reflection Equation} algebra
${\cal L}(R)$
is generated by the set of $N^2$ generators $L_i^{\, j}$
satisfying the condition
\be
\lb{RLRL}
R_{12} L_1 R_{12} L_1 = L_1 R_{12} L_1 R_{12} \; .
\ee
This algebra is naturally endowed with the structure of adjoint comodule
under the left coaction of quantum linear group, and with the
structure of trivial (invariant) comodule under its right coaction:
\be
\lb{coac}
\delta_{\ell}(L_i^{\, j}) = T_i^{\, k}
{\cal S}(T_p^{\, j}) \otimes L_k^{\, p} \;\;\; ,
\;\;\; \delta_r(L_i^{\, j}) = L_i^{\, j} \otimes 1 \; .
\ee
It is worth noticing that
in case of $R$ being the ${\cal R}$-matrix
corresponding to a linear quantum group,
the matrix elements $L_i^{\, j}$ have a nice geometrical
interpretation as quantum analogues of the basic
right-invariant vector fields over the linear group.

\noindent
{\bf Quantum trace.}

Let us derive several consequences of the relation (\ref{RTT}).
Denote ${\cal R}_{12} \equiv P_{12} R_{12}$, where $P_{12}$ is the
permutation matrix, and rewrite (\ref{RTT}) in a more conventional
form
\be
\lb{RTT2}
{\cal R}_{12} T_1 T_2 = T_2 T_1 {\cal R}_{12} \; .
\ee
We say that  the Hecke symmetry is {\em closed} if the matrix
${\cal R}^{t_1}_{\, 12}$ is invertible. With this additional assumption
one can turn (\ref{RTT2}) to the form
\be
\lb{RTT3}
({{\cal R}^{t_1})^{-1}}_{i_1 i_2}^{\, k_1 k_2}
{{\cal S}(T)}_{k_2}^{\, j_2} {T^t}_{k_1}^{\, j_1} =
{T^t}_{i_1}^{\, k_1} {{\cal S}(T)}_{i_2}^{\, k_2}
({{\cal R}^{t_1})^{-1}}_{k_1 k_2}^{\, j_1 j_2} \; .
\ee
Now, by performing the antipodal map of the above relation
and summing up the indices $j_1$, $j_2$ we eventually
obtain
\be
\lb{square1}
{\cal S}^2(T) {\cal C} = {\cal C} T \; , \quad \mbox{where} \quad
{\cal C}_i^{\, j} \equiv \sum_k ({{\cal R}^{t_1})^{-1}}_{j i}^{\, k k} =
Tr_{(1)}
\left[ \left( ({\cal R}^{t_1})^{-1} \right)^{t_1}_{\, 12} P_{12} \right]
\ee
An analogous procedure with summing up indices $i_1$, $i_2$ gives
\be
\lb{square2}
{\cal B} {\cal S}^2(T) = T {\cal B} \; , \quad \mbox{where} \quad
{\cal B}_i^{\, j} \equiv \sum_k ({{\cal R}^{t_1})^{-1}}_{k k}^{\, i j} =
Tr_{(2)}
\left[ \left( ({\cal R}^{t_1})^{-1} \right)^{t_1}_{\, 12} P_{12} \right]
\ee
By definition, all the entries of the matrix $T$ are linearly independent
and, hence, the matrix ${\cal BC}={\cal CB}$ is scalar.

With the use of the matrix $\cal C$ one introduces an analogue of the
trace operation in the algebra ${\cal L}(R)$ \cite{Resh}.
Namely, the {\em quantum trace} of the matrix $L$ is defined as
\be
\lb{tr}
Tr_q L \equiv Tr ({\cal C} L)\; ,
\ee
and it extracts the invariant
part of the left-adjoint comodule $L_i^{\, j}$.
Indeed, using formula (\ref{square1}) one easily checks
\ba
\nonumber
{\delta}_{\ell}(Tr_q L) &=& Tr \left\{ {\cal C} T L {\cal S}(T) \right\} =
Tr \left\{ {\cal S}^2(T) {\cal C} L {\cal S}(T) \right\} \\
\lb{scal}
&=&
{\cal S}\left( Tr \left\{ T {\cal S}(T) {\cal C} L \right\} \right) =
\id \otimes Tr_q L \; .
\ea
Here we omit symbol $\otimes$ in calculations. The above reasonings
also work
for any left-adjoint comodule, for example, the quantum trace
of any power of the matrix $L$ is invariant.
%Analogously, the matrix $\cal B$ can be used to extract the scalar
%part of the right-adjoint comodules, like ${\cal S}(T)LT$.

When realizing (\ref{square1}) and (\ref{scal}) in specific
$T$-representations, one gets several useful relations for
the quantum trace and the matrix $\cal C$.
Consider the representations
$\rho(T_a) = {\cal R}_{ab}$, and $\rho'(T_b) = {\cal R}^{-1}_{ab}$.
Then
$$
\begin{array}{lll}
\rho({\cal S}(T_a)) =
\left( \left( ({\cal R}^{-1}_{ab})^{t_a}\right)^{-1} \right)^{t_a}
& , &
\rho({\cal S}^2(T_a)) = {\cal R}^{-1}_{ab} \; , \\
\rho'({\cal S}(T_b)) = {\cal R}_{ab} & , &
\rho'({\cal S}^2(T_b)) =
\left( ({\cal R}^{t_b}_{ab})^{-1} \right)^{t_b} =
\left( ({\cal R}^{t_a}_{ab})^{-1} \right)^{t_a} \; .
\end{array}
$$
The relation (\ref{square1}) in the representation $\rho$ looks like
$$
{\cal C}_a {\cal R}_{ab} =
\left( \left( ({\cal R}^{-1}_{ab})^{t_a}\right)^{-1} \right)^{t_a}
{\cal C}_a \; ,
$$
whereas for the representation $\rho'$ it can be cast into the form
$$
{\cal R}_{ab} {\cal C}_b = {\cal C}_b
\left( \left( ({\cal R}^{-1}_{ab})^{t_a}\right)^{-1} \right)^{t_a} \; .
$$
Combining these two relations we obtain
\be
\lb{RCC}
R_{12} {\cal C}_1  {\cal C}_2 = {\cal C}_1 {\cal C}_2 R_{12} \; .
\ee
In a similar way one deduces from (\ref{scal})
\ba
\lb{scalR}
&& {Tr_q}_{(2)} ( R_{12} X_1 R_{12}^{\, -1} ) =
{Tr_q}_{(2)} ( R_{12}^{\, -1} X_1 R_{12} ) = (Tr_q X) \id \; , \\
\lb{scalR2}
&& {Tr_q}_{(1,2)} ( R_{12} X_{12} R_{12}^{\, -1}) =
{Tr_q}_{(12)} X_{12} \; .
\ea
Here $X_1$ and $X_{12}$  are arbitrary operator-valued matrices, and
the last relation is a consequence of (\ref{RCC}). Now, applying
({\ref{scalR}) to the relations (\ref{RLRL}) one easily checks
that all the elements $Tr_q L^k$, $k\geq 1$ are central
in the algebra ${\cal L}(R)$.
\bigskip

\noindent
{\bf Quantum antisymmetrizers and quantum Levi-Civita tensors.}

Now let us remind some of the results of the paper \cite{Gur} which
will be relevant to considerations below.
Note that the notation adopted here is slightly different from that of
Ref.\cite{Gur}.
The correspondence is the following:
our parameter $q$ corresponds to $q^{1/2}$ of Ref.\cite{Gur};
$R$ corresponds to ${1 \over q^{1/2}} S$ of Ref.\cite{Gur}.

By its definition any Hecke symmetry determines a
series of local representations of Hecke algebras $H_q(\C)$
(see e.g. \cite{Wenzl}).
We will consider some central projectors of these
Hecke algebras, namely, antisymmetrizers ${P_-}^k$.
They are iteratively defined by the relations
$$
{P_-}^1 = \id \; , \quad {P_-}^k = {1 \over k_q} \left( q^{k-1} \id
-q^{k-2}R_{k-1} + \dots + (-1)^{k-1} R_1 \cdot \dots \cdot R_{k-1} \right)
{P_-}^{k-1} \; .
$$
Here we use the notation $R_i \equiv R_{i (i+1)}$ and ${P_-}^i \equiv
({P_-}^i)_{12\dots i}$. One can give several alternative definitions
for the antisymmetrizers:
$$
{P_-}^k = {1 \over k_q} {P_-}^{k-1} \left( q^{k-1} \id - q^{k-2}R_{k-1}
+ \dots + (-1)^{k-1} R_{k-1} \cdot \dots \cdot R_1 \right) \; ,
$$
or
$$
{P_-}^k = {1 \over k_q} \left( q^{k-1} \id
-q^{k-2}R_1 + \dots + (-1)^{k-1} R_{k-1} \cdot \dots \cdot R_1 \right)
{P_-}_2^{k-1} \; ,
$$
where ${P_-}^i_j \equiv ({P_-}^i)_{j(j+1)\dots (j+i-1)}$.

The following properties of the
projectors ${P_-}^k$ are proved in \cite{Gur}:
\ba
\lb{1}
{P_-}^k R_i &=& R_i {P_-}^k = - {1 \over q_{_{_{}}} } {P_-}^k \; ,
\quad \mbox{for} \quad i \leq k-1 \; ;
\\
\lb{2}
{P_-}^k {P_-}_j^i &=& {P_-}_j^i {P_-}^k = {P_-}^k \; ,
\quad \mbox{for} \quad i+j-1 \leq k \; ;
\\
\lb{3}
{P_-}^k R_k {P_-}^k &=& - { (k+1)_q^{^{^{}}} \over k_{q_{_{}}} }
{P_-}^{k+1} + { q^k \over k_q} {P_-}^k \; ;
\\
\lb{4}
{P_-}_2^k R_1 {P_-}_2^k &=& - { (k+1)_q^{^{^{}}} \over k_{q_{_{}}} }
{P_-}^{k+1} + { q^k \over k_q} {P_-}_2^k \; .
\ea
The last pair of relations may be equally used for
an iterative definition of ${P_-}^k$.
\bigskip

The closed Hecke symmetry $R$ is called the {\em even Hecke symmetry
of rank p} if the following condition on antisymmetrizer
is satisfied
%\footnote{
%For motivation and discussion of this definition see \cite{Gur}.
%}
\be
\lb{5}
dim {P_-}^{p+1} = 0 \; .
\ee
As a consequence of (\ref{5}) one gets for the closed symmetry
that the image of ${P_-}^p$ is one-dimentional \cite{Gur}:
$ dim {P_-}^p = 1 $.
Therefore ${P_-}^p$ may be presented in the form
\be
\lb{epsilon1}
({P_-}^p)_{i_1 i_2 \dots i_p}^{\, j_1 j_2 \dots j_p} =
u_{i_1 i_2 \dots i_p} v^{j_1 j_2 \dots j_p} \equiv \u \vv \; ,
\ee
where
$$
v^{i_1 i_2 \dots i_p} u_{i_1 i_2 \dots i_p} = 1 = \vv \u \; ,
$$
and in the r.h.s. of the formulas we introduce the symbolic notation
to be employed in what follows. The quantities $u_{|\dots{\cal i}}$ and
$v^{{\cal h} \dots |}$ are the analogues of left and right
Levi-Civita tensors.
They are defined (up to a numerical factor) by the conditions
\be
\lb{epsilon2}
R_i \u = - {1 \over q} \u \; , \qquad \vv R_i = - {1 \over q} \vv \; ,
\qquad  i = 1, \dots  ,p-1 \; .
\ee
One can easily check that the tensors $u_{|\dots{\cal i}}$,
$v^{{\cal h} \dots |}$ defined by (\ref{epsilon2})
really enter into formula (\ref{epsilon1}) for ${P_-}^p$ (simply
apply ${P_-}^p$ to these $u_{|\dots{\cal i}}$ and
$v^{{\cal h} \dots |}$).

As a consequence of (\ref{3}), (\ref{4}) and the definitions of
$\cal C$ and $\cal B$ (\ref{square1}), (\ref{square2}) we have
\be
\lb{BC}
{\cal C}_{|a{\cal i}}^{\; \; \, {\cal h}b|} = { p_q \over q^p}
v^{{\cal h} b 1 \dots p-1 |} u_{| a 1 \dots p-1 {\cal i}} \; , \qquad
{\cal B}_{|a{\cal i}}^{\; \; \, {\cal h}b|} = { p_q \over q^p}
v^{{\cal h}  1 \dots (p-1) b|} u_{|  1 \dots (p-1) a{\cal i}} \; ,
\ee
and also
\be
\lb{trR}
Tr_q \id = Tr {\cal C} = { p_q \over q^p } = Tr {\cal B} \; , \qquad
{Tr_q}_{(2)} R_{12} = Tr_{(1)} {\cal B}_1 R_{12} = \id \; .
\ee

Let us derive a few more useful relations between the quantum trace and the
quantum Levi-Civita tensors. For any (possibly operator-valued) $N\times N$
matrix $X_i^{\, j}$ define the {\em symmetrizing map} ${S_+}^k$:
$Mat(N) \rightarrow Mat(N)^{\otimes k}$
\be
\lb{S}
{S_+}^k(X) \equiv X_1 + R_1 X_1 R_1 + \dots + R_{k-1} R_{k-2} \cdot \dots
\cdot R_1 X_1 R_1 R_2 \cdot \dots \cdot R_{k-1} \; .
\ee
Obviously, for any $X$,
$$
{S_+}^k(X)  R_i =  R_i {S_+}^k(X)\; ,
\quad i = 1, \dots ,k-1 \; .
$$
Now performing simple calculations one obtains the formula
$$
{P_-}^p {S_+}^p(X) = {S_+}^p(X) {P_-}^p = { p_q \over q^{p-1}}
{P_-}^p X_1 {P_-}^p \; ,
$$
whereof the desired relations follows:
\be
\lb{useful}
\vv {S_+}^p(X) = q Tr_q X \, \vv  \; , \qquad
{S_+}^p(X) \u = \u \, q Tr_q X \; .
\ee
\bigskip

\section{Newton relations and Cayley-Hamilton theorem}
\setcounter{equation}0

For any RE algebra
associated with rank $p$ even Hecke symmetry
there are two canonical possibilities for introducing
a set of $p$ central elements being at the same time
invariants of the adjoint coaction $\delta_{\ell}$ (\ref{coac}).
The first set contains the $q$-traces of powers of the matrix $L$
\be
\lb{s}
s_q(i) \equiv q \, Tr_q L^i \;, \quad i= 1, \dots ,p.
\ee
The centrality and invariance of these elements have been already
demonstrated in the previous Section.

The second set is not so obvious as the first one.
It is formed by combinations of quantum
minors of the matrix $L$:
\be
\lb{sigma}
\sigma_q(i) \equiv \alpha_i \vv
\left( L_1 R_1 \dots R_{i-1} \right)^i \u \; .
\ee
Here $\alpha_i$ are some normalizing constants to be fixed below.
As a partial explanation to this definition
one may note that the combinations of ${\cal R}$-matrices
enter the r.h.s. of
(\ref{sigma}) namely in such a way that the whole expression
would be adjoint invariant. One can also directly check the
centrality of the elements $\sigma_q(i)$, but
this will follow from the connection between the first and the
second sets established below.

The elements $s_q(i)$ and $\sigma_q(i)$ play,
respectively,  the role
of the basic power sums and the basic symmetric polynomials
for the quantum matrix $L$.
In case of the quantum linear groups and
for $q$ being not the root of unity
both sets are known to be complete, i.e.,
they generate the center of RE algebra.
As for other Hecke symmetries,
it seems that the problem of completeness of
these sets should be treated in each case separately.
The relation between these two sets is established by the
following $q$-analogue of Newton's formulas:

{\bf Proposition 1.}\hspace*{5 true mm}
\em
The sets  $\{ \sigma_q(i) \}$ and $\{ s_q(i) \}$,  $i=1, \dots ,p$
are connected by the relations
\be
\lb{newton}
{i_q \over q^{i-1}} \sigma_q(i) -
s_q(1)\,\sigma_q(i-1) + \dots +
(-1)^{i-1} s_q(i-1) \,\sigma_q(1) +
(-1)^i s_q(i) = 0 \, ,
\ee
provided that the numerical factors $\alpha_i$ are fixed as
\be
\lb{alpha}
\alpha_i = q^{-i(p-i)} { p \choose i}_q  \; .
\ee
Here ${ p \choose i}_q = i_q! (p-i)_q! / p_q!$
--- are the $q$-binomial coefficients, and
$i_q! = i_q \, (i-1)_q!$ ,  $1_q! = 1$.
\rm
\par
{\bf Proof:}
Consider the quantities $s_q(i-j)\sigma_q(j)$ for
$1 \leq j < i \leq p$.
With the help of (\ref{useful}) and
taking into account the commutativity relations
$$
R_k ( L_1 R_1 \dots R_{j-1} )^j = ( L_1 R_1 \dots R_{j-1} )^j R_k \; ,
\quad k=1, \dots ,j-1 \; ,
$$
following from (\ref{YB}), (\ref{RLRL}),
one can perform transformations:
\begin{eqnarray*}
s_q(i-j) \sigma_q(j) &=&
\alpha_j \vv  {S_+}^p(L^{i-j}) (L_1 R_1 \dots R_{j-1})^j \u \; = \\
&&
\alpha_j {j_q \over q^{j-1}} \, \vv (L_1^{i-j+1} R_1\dots R_{j-1})
(L_1 R_1 \dots R_{j-1})^{j-1} \u \; + \\
&&
\alpha_j {(p-j)_q \over q^{p-j-1}} \, \vv
(R_j \dots R_1 L_1^{i-j}R_1\dots R_j)
(L_1 R_1 \dots R_{j-1})^j \u  \; .
\end{eqnarray*}
Now using relations
$$
(L_1 R_1 \dots R_{j-1})^j R_j \dots R_1 = (L_1 R_1 \dots R_j)^j
$$
we complete the transformation of the second term of the sum,
and finally have
\begin{eqnarray}
\nonumber
s_q(i-j) \sigma_q(j) &=&
\alpha_j {j_q \over q^{j-1}} \, \vv (L_1^{i-j+1} R_1\dots R_{j-1})
(L_1 R_1 \dots R_{j-1})^{j-1} \u \; + \\
\lb{s-sigma}
&&
\alpha_j {(p-j)_q \over q^{p-j-1}} \, \vv
(L_1^{i-j}R_1\dots R_j)
(L_1 R_1 \dots R_j)^j \u  \; .
\end{eqnarray}
In the boundary cases $j=1$, $j=i-1$ this equation reads
\begin{eqnarray}
\nonumber
s_q(i-1) \sigma_q(1) &=&
\alpha_1 { q^{p-1} \over p_q} s_q(i) +
\alpha_1 {(p-1)_q \over q^{p-2}} \, \vv
(L_1^{i-1}R_1)
(L_1 R_1) \u  \; ; \\
\nonumber
s_q(1) \sigma_q(i-1) &=&
\alpha_{i-1} {(i-1)_q \over q^{i-2}} \, \vv (L_1^{\, 2} R_1\dots R_{i-2})
(L_1 R_1 \dots R_{i-2})^{i-2} \u + \\
\nonumber
&&
{\alpha_{i-1} \over \alpha_i} {(p-i+1)_q \over q^{p-i}} \sigma_q(i) \; .
\end{eqnarray}
Let us choose $\alpha_1 = p_q/q^{p-1}$ in order that
$\sigma_q(1) = s_q(1)$
Now putting
$$
\alpha_j = q^{2j-1-p} { (p-j+1)_q \over j_q } \alpha_{j-1}\; ,
$$
and, then, making an alternating sum of (\ref{s-sigma})
for different values of $1 \leq j \leq i-1$ one obtains
the relation (\ref{newton}).
~~\rule{2.5mm}{2.5mm}\par

The elements $s_q(k)$ for $k > p$ can also be expressed
in terms of $\{ \sigma_q(i) \}$, or $\{ s_q(i)\}$, $ i \leq p$.
The relation is provided by an analogue of the Cayley-Hamilton
Theorem:

{\bf Proposition 2.}\hspace*{5mm}
\em
The $L$-matrix satisfies the following characteristic identity
\be
\lb{char}
\sum^{p}_{i=0} (-L)^i \sigma_q(p-i) \equiv 0 \; .
\ee
Here we imply $\sigma_q(0) = 1$.
\rm
\par
{\bf Proof:}
Consider the quantity
$$
w(x)^{{\cal h}1 2 \dots p|} \equiv  \prod_{i=0}^{p-1}
\left[ (L_1 -q^{2i}x \id ) R_1 \dots  R_{p-1} \right] \u \; .
$$
Generalizing in an obvious way an observation
that the commutator
\par
\noindent
$ \left[ R_1 \, , \, (L_1 -x\id) R_1 (L_1 - q^2 x \id) \right] $
is proportional to the $q$-symmetric term $(R_1 + 1/q)$,
one can convince himself that
$R_i w(x)^{{\cal h}1 2 \dots p|} = 0$ for $i = 1, \dots,p-1$.
Hence
\be
\lb{w}
w(x)^{{\cal h}1 2 \dots p|} =
{\bigtriangleup}(x) \vv \; .
\ee
Here the scalar coefficient
${\bigtriangleup}(x) \equiv w(x)^{{\cal h}1 2 \dots p|} \u$
is an analogue of the characteristic polynomial for the matrix $L$.
It can be calculated with the use of (\ref{YB}), (\ref{RLRL}),
(\ref{epsilon2}), and the $q$-combinatorial relations. The result is
\be
\lb{polynom}
{\bigtriangleup}(x) =
\sum^{p}_{i=0} (-x)^i \sigma_q(p-i) \; .
\ee
Now, simply repeating the classical arguments of \cite{Lank}
one can prove that ${\bigtriangleup}(L) \equiv 0$.
The starting point for this is provided, e.g., by the relation
$$
w(x)^{{\cal h}2 \dots p+1|} \u =
{\bigtriangleup}(x)
v^{{\cal h}2 \dots p+1|} u_{|1 \dots p{\cal i}} \; .
$$
~~\rule{2.5mm}{2.5mm}\par
To conclude it is worth making one note. In the case
$q\rightarrow 1$ (that is $R^2=1$)  formulae (\ref{newton})
and (\ref{char}) tend to their classical limits. We would
like to emphasize that the key role is played by the value
of {\it rank p\/} of  Hecke  symmetry. The specific value $N$
of the dimension of the space $V$ where this symmetry
is realized does not matter here.

\section*{Acknowledegments}
This work was initiated when two of the authors (D.G. and  P.P.)
participated at the Stefan Banach Center Minisemester on
{\em 'Quantum Groups and Quantum Spaces'} (Autumn 1995, Warsaw).
Its a pleasure to acknowledge the warm hospitality of the
organizers of the Minisemester.
%and especially to S.Zakrzewski and S.L.Woronowicz
%for providing us with warm and creative atmosphere.
We are indebted to A.Isaev and G.Arutyunov for numerous
valuable discussions.
The work of P.P. was supported in part by the RFBR grant 95-02-05679a
and by the INTAS grant 93-2492 ( the research program of the ICFPM).

\end{document}